# Universality Properties of Gaussian Quadrature, The Derivative Rule, and a Novel Approach to Stieltjes Inversion

William P. Reinhardt[*]
Department of Chemistry[**]
University of Washington, Campus Box 351700
Seattle, WA 98195-1700

**Abstract:** In the 1890s Stieltjes, Perron, and Markov in analysis of the moment problem indicated that a non-decreasing measure, $\mu(x)$, consistent with the moments and absolutely continuous on a real interval [a,b], could be constructed as the discontinuity across branch cut of a Cauchy representation of an otherwise real-analytic function $F(z) = \int_a^b d\mu(x)/(z-x)$. $F(z)$ being constructed as a continued fraction from the traditional three-term recursion relation related to the polynomials, $P_n(x)$, orthogonal on [a,b] with respect to $d\mu(x)$, leading to the inversion process of Markov, Perron and Stieltjes. Introduction of Gaussian quadratures suggests more direct approaches, an idea which is culminated herein. The real, positive, and continuous weight function, $\rho(x) = \mu'(x)$ is constructed directly from interpolation of the zeros, $x[k]$, of the $P_n(x)$, ordered as an increasing function of the integral values, k, and the corresponding Gaussian quadrature weights, $w[k]$, often referred to as the Christoffel numbers. The success of this novel numerical approach rests on a simple theorem: namely, that the desired weight function, $\rho(x[k])$, evaluated at the quadrature point $x[k]$, is approximated by $w[k]/(dx[k]/dk)$, with "k" now being assumed to be a continuous real variable allowing said derivative to be calculated from a smooth interpolation of the quadrature abscissas at their integer values. This *derivative rule* is shown to be closely related to a simple type of *universal behavior*, relating the *clock-rule* spacings of the $x[k]$ to the derivatives $dx[k]/dk$ which, once established, allows an immediate derivation of the derivative rule, via known and well understood *universal* clock-rule results for Nevai-Blumenthal (N-B) class polynomials. This framework also leads to a statement *of a universality for quadrature weights and weight functions*. This is implicit in statements of clock-rule universality for the N-B class, but not usually made explicit, as the normal utility of the clock rule is not Stieltjes inversion, but the understanding of the distributions of zeros of random polynomials. This new universality is illustrated for several classical orthogonal polynomials, and the non-classical Pollaczek polynomials. The derivative rule is then applied to the problem of numerical Stieltjes inversion on both finite and unbounded open intervals, and found to give exponential convergence in *n* for the weight function. Exponential convergence stands in contrast to the power law convergence of historically well-known histogram techniques which may be traced back to Markov, Chebyshev, and Stieltjes. A final overview Section traces the origins of the derivative rule and its use in developing novel numerical approaches to the use of approximate spectral resolutions of Schrödinger operators found using $L^2$ basis-sets, illustrated here in computation of a photo-effect cross section.

[*]Professor Emeritus.
[**]Email: rein@chem.washington.edu

# 1. Introduction: Overview and Introductory Material

## A) An Overview

In this Introduction we review the basic ideas of numerical solution of the classical *Problem of Moments*, going back to the days of Stieltjes, Markov, Chebyshev, and Perron in the late 19th Century, followed by a discussion of Gaussian quadrature, equivalent quadrature, and the relationship Gaussian quadrature to the Markov continued fraction. Section 2 begins with an outline of notations and assumptions. The *derivative rule conjecture* is introduced and indications of its utility for Stieltjes inversion, and evaluation of the *equivalent quadrature weights* of Section 1 presented. The derivative rules is then shown to be exact for the Chebyshev polynomials where analytic expressions for the quadrature abscissas and weighs are available, as summarized in the Appendix. Section 3 posits universal relationships leading to a new approach to Stieltjes inversion. This follows from interpreting universal results, such as clock-rule behavior of polynomial zeros, in a simple but novel manner, as derivatives. The resulting *derivative rule* relates the weight function generating the orthogonal polynomials directly to the zeros of those polynomials and the Gaussian quadrature weights, or Christoffel numbers, associated with them. A universality involving the integration weights and the corresponding weight functions then follows and is illustrated.

Section 4 A reviews the traditional numerical histogram Stieltjes inversion methodology where convergence is typically $\sim 1/N^2$. Section 4B approaches the same problems via the derivative rule, and seen, to give convergence as approximately $10^{-N}$, for Chebyshev systems, using only the Gaussian quadrature weights and abscissas, as input This is extended, in Section 4C), to the determination of a Gegenbauer weight function from the now fully numerically determined Gaussian weights and abscissas, and where exponential convergence of the form $10^{-N/\alpha}$, with $\alpha \sim 3$ is found. This example provides an illustration of problems of the required interpolations as used near boundaries, and an initial practical solution indicated. In Section 4D the restriction to a bounded closed interval is removed, and then illustrated for determination of the Hermite weight function where exponential convergence of the form $10^{-N/\alpha}$, also with $\alpha \sim 3$, is found.

The origins of these ideas are overviewed in Section 5, as they arose from solution of strikingly similar mathematical problems in utilization of discretized spectral resolutions of Schrödinger operators. These ideas having been introduced as useful computational approaches in early and mid 1970s. A summary and conclusions appear in Section 6.

## B) Orthogonal Polynomials and the Moment Problem

In Sections 1-4, the focus is on the case of a non-decreasing, absolutely continuous, real function μ(x) on a single compact real interval x ∈ [a,b], although in one occasion μ(x) may have jump discontinuities outside of [a,b], in which case the integration range would be extended. It is assumed that all moments,

$$\mu_n = \int_a^b x^n d\mu(x) \tag{1}$$

are real and finite for all non-negative n.  Orthogonal polynomials $P_n(x)$, often referred to as OPs, being real linear combinations of the $x^n$ satisfy

$$\int_a^b P_n(x) P_m(x) d\mu(x) = 0, \quad n \neq m, \quad n,m = 0,1,2... \tag{2}$$

and are typically standardized via the canonical three term recursion, Ref.(1),

$$P_{n+1}(x) = (A_n x + B_n) P_n(x) - C_n P_{n-1}(x), \quad P_0 = 1, \quad P_1 = A_0 x + B_0. \tag{3}$$

The recursion coefficients satisfy the positivity product condition $(A_{n-1} A_n C_n) > 0$ for $n > 0$, as required by Favard's theorem.  This latter being a sufficient condition for use of Eqn. (3) and implying the existence of a positive measure such that the polynomials, as so determined, are orthogonal. This is the standard foundations of the theory of orthogonal polynomials; see Refs. (1-3).  The set of definitions and relationships in Eqns. (1-3) will be referred to an OP system.

The moment inversion problem, which is assumed to be well determined, (see Refs. (3,4)), may then be stated as:  Given the moments of Eqn. (1), or the recursion coefficients of Eqn. (3), can we determine, numerically, if not analytically, the Stieltjes measure $d\mu(x)$ on [a,b]?  Or, as in most of what follows, can we determine the non-negative weight function, $\rho(x)$, $\rho(x) = \mu'[x]$, in the case that $\mu(x)$ is absolutely continuous on [a,b]?

The traditional inversion process of Markov and Stieltjes is an exercise in analytic continuation combined with continued fractions, Refs. (4-7), determined by the recursion coefficients of Eqn. (3),

$$F(z) = \lim_{n \to \infty} A_0 / (A_0 z + B_0 - C_1/(A_1 z + B_1 - C_2/(A_2 z + B_2 \ldots - C_n/(A_n z + B_n)), \tag{4}$$

showing that, in that limit, F(z) converges to a real analytic function of form, for $z \notin [a,b]$,

$$F(z) = \int_a^b d\mu(x)/(z-x). \tag{5}$$

When truncated at the $n^{th}$ convergent, n being the last denominator in Eqn. (4), the resulting truncation, $F_n(z)$, of the continued fraction, Eqn. (4) can be written as a ratio degree n-1 and n polynomials:

$$F_n(z) = P_n(z,1) / P_n(z), \tag{6}$$

Where the $P_n(z)$ are the polynomials of Eqn. (3), and the $P_n(z,1)$ are the related numerator OPs, see Ismail Ref. (2), sections 2.3 and 2.6. Eqn. (6) implies that the zeros of the $P_n(x)$ are poles of the $n^{th}$ approximant, a result famous for its role in the early developments of the theory of orthogonal polynomials, see, again, Refs. (1-3). Note also that as the zeros of the $P_n(x)$ all lie in the interval [a,b], that any finite $n$ approximation, via Eqns. (5) or (6), will have $n$ poles in that compact interval where, as $n \to \infty$, these merge to form a (branch cut) discontinuity along [a,b] in the limiting integral. This property of such an approximation is well displayed in Eqn. (10) in the following Subsection.

For the case where $d\mu(x) = \rho(x)dx$, and $\rho(x)$ is a positive continuous real function on [a,b], a simplified form of the Markov-Stieltjes inversion formula, useful in numerical work, then is

$$\lim_{z \to x' + i\varepsilon} \int_a^b d\mu(x)/(z-x) \;=\; \lim_{z \to x' + i\varepsilon} \int_a^b \rho(x)\,dx/(z-x) \tag{7a}$$

$$= \mathcal{P}\int_a^b \rho(x)\,dx/(x'-x) \;-\; i\pi\rho(x') \tag{7b}$$

where $\varepsilon$ is a positive infinitesimal, vanishing in the limit, $x' \in [a,b]$, and $\mathcal{P}$ denotes the Cauchy principal value of the integral.

It is the second form of Eqn. (7a,b) which is commonly used in physics applications as arise, for example, in quantum scattering theory, see Refs. (8,9) and which will be used making a connection between methods of use in the classical moment inversion problem to the closely related spectral theory of Schrödinger operators in Section 5.

Numerical approximations, using an approximate $F_n(z)$, may converge well away from the cut on [a,b], and then be numerically analytically continued back to the real axis. Such numerical continuation, say by use of a lower order continued fraction, or a Padé approximant, Refs.(10 - 14), which when fit to data away from the real axis, nonetheless, when evaluated in the $\pm i\varepsilon$ limits, gives quite useful numerical results. Gautschi, Ref. (15), discusses an explicit discussion of a numerical $\pm i\varepsilon$ limit, using an accelerated convergence technique, see Temme, Ref. (16). In what follows, unless otherwise indicated, it is assumed that only an absolutely continuous weight function is considered.

**B) Gaussian Quadrature and Equivalent Quadratures**

    **i) Gaussian Quadrature**
Everything to follow depends on knowledge of the method of *Gaussian mechanical quadrature*, Refs. (1,2,3), which will simply be referred to as *Gaussian quadrature*. Gaussian quadrature, of order n, requires knowledge of $P_n(x)$, the orthogonal polynomials generating by *weight function* $\rho(x)$, the *n quadrature abscissas* x[k], or x[k,n], which are the zeros of the $n^{th}$ degree OPs, and the *n quadrature weights* w[k], or w[k,n], also referred to as the Christoffel numbers, associated with these zeros.  Note that "*n*" as part of the notation for the x[k], and w[k] will only be included where ambiguity would otherwise arise, as, for example in the important case where the limiting behavior as $n \to \infty$ is under discussion.

The usual *Gaussian quadrature* procedure is as follows.  An integral of the type

$$\int_a^b g(x)\rho(x)dx. \tag{8}$$

is approximated as

$$\int_a^b g(x)\rho(x)dx = \sum_{k=1}^n w[k]\, g(x[k]) \tag{9}$$

This n-point approximation gives an exact result if g(x) is, itself, a polynomial in x of degree 2n-1 or less.   If g(x) is not a polynomial the result will not be exact, but may be well approximated to the extent that g(x) is well approximated by a linear superposition of the polynomials $P_m(x)$ for m ≤ 2n-1.

### ii) Equivalent Quadrature

Suppose the *Gaussian quadrature* approximation of Eqns. (8,9) is repeated to perform that quadrature, with the same Gaussian abscissas, x[i], and weights, w[i] and weight function ρ(x), all on the compact interval [a,b], to approximate the integral

$$\int_a^b g(x)dx, \tag{10}$$

where now ρ(x) is conspicuously absent.  Carrying out this process, *indeed mechanically,* is straightforward: as now the (non-negative) ρ(x) of Eqn. (15c) is absent, but still wishing to use the Gauss quadrature abscissas and weights generated by ρ(x), the integral is approximated as

$$\int_a^b g(x)dx = \int_a^b \rho(x)\,(g(x)/\rho(x))dx = = \sum_{i=1}^n w[i]\, g(x[i])/\rho(x[i])$$

$$= \sum_{i=1}^n w_i^{EQ}\, g(x[i]), \tag{11}$$

where the *Equivalent Quadrature weight* is defined as:

$$w_k^{EQ} = w[k]/\rho(x[k]). \tag{12}$$

The definition of Eqn. (12) follows automatically from the formulaic application of the rules of Gaussian quadrature. Nothing has been proven here, standard formulae just rearranged. Why is this of use? First of all, application of Eqn. (11) often gives excellent numerical results, even though only being exact if g(x)/ρ(x) is a polynomial of degree 2n-1 or less.

However, the real utility of Eqn. (11) arises when (as in the quantum mechanical example discussed in Section 5) a Gaussian quadrature arises naturally, but not one related to having its weight function present in the integral to be performed. What will then, perhaps surprisingly, be seen is that $w_k^{EQ}$, via the yet to be introduced *derivative rule*, can often be well approximated without specific knowledge of either the w[k] or ρ(x[k]).

### iii) The Markov Approximation as a Gaussian Quadrature

The approximation of replacing F(z) by the $F_n(z)$ of Eqn. (6) may be re-written as (see Chihara, Ref (2), pages 88-89, Shohat and Tamarkin, Ref(4), pp 39-43 , and Ismail, Ref(2), Sections 2.5 and 2.6) in terms of the now introduced quadrature abscissas and weights generated by ρ(x):

$$F(z) = \int_a^b \rho(x)dx/(z-x) \cong F_n(z) = \sum_{k=1}^n \frac{w[k]}{z-x[k]} \ . \quad (13)$$

This simplifies understanding of the representations of Eqns. (4 and 6) as approximations to the integral, as well as leading quickly to the traditional histogram Stieltjes inversion methods outlined in Section 4A), which may be traced back to Chebyshev, see the above references, but where the simpler pedagogical development of Chihara is followed.

### 2) Notations, The Derivative Rule, Simple Examples

Notations and the derivative rule are introduced. The possible utility of the derivative rule in Stieltjes inversion is indicated. The derivative rule is shown to be exact for the Chebyshev families of polynomials, setting the stage for the derivations which appear in Section 3. Descriptions of the Chebyshev OP systems appear in the Appendix.

#### A) Notations

i) In what follows we will often consider x[k], and sometimes even w[k], to be continuously differentiable functions of the "quadrature index k," so these are not labelled in the conventional manner as $x_k$ and $w_k$. Unless otherwise noted, x[k] ∈ [a,b].

ii) When the dependence of x[k] or w[k] on *n* is of importance, as in statements *of universality* in Section 3, these will be denoted as x[k,n] and w[k,n] respectively.

iii) The notation of the non-negative continuous weight function, ρ(x), on the single interval, x∈ [a,b], is used rather than the conventional w(x), as for example in Refs. (1a,b). This is to avoid possible confusion with the quadrature weights, w[k], as labelled by integer values of k.

iv) The derivative of x[k] for , x[k] ∈ [a,b], is taken in the sense defined by

$$x'[k] = (dx[i]/di)_{i=k} \qquad (14)$$

where *i* is a continuous variable and which indicates either an analytic or numerically smooth interpolation of the discrete set x[k].

v) The x[k,n] are ordered as increasing functions of *k*, for fixed *n*, giving a positive derivative x'[k,n], and the x'[k,n] then, most often, only used for integer values of k. Note that an opposite ordering often appears in statements of *universality* and the *clock-rule*, as discussed in Section 3, below. An exception to the use of only integer *k* is Figures 1 and 2, where both x[k,n] and x'[k,n] are known as a continuous functions of *k*, for fixed *n*.

vi) In Sections 4 and 5 which focus on numerical work: N, rather than the n of Sections 1-3, will be used to facilitate readability of Tables and Figures.

## B) The Derivative Rule and Stieltjes Inversion

The *derivative rule* is easy to express, see Eqn. (14) for a comment on the notation,

$$x'[k] = w[k]/\rho(x[k]) = w_k^{EQ} , \qquad k = 1,2,...n. \qquad (15a)$$

This assertion was introduced in the early 1970s, by Heller, Ref. (17), and used by Yamani and Reinhardt, Ref. (18, Section V.B) who refer to this as *Heller's derivative rule*, for calculation of the above *equivalent quadrature weights* of Eqn. (12) in the context of constructing the real and positive discontinuity in an integral such as in Eqn. (7), as it might appear in the spectral resolution of a *Schrödinger operator*, see Refs (19, 20). Details of this early work, with further references, appear in Section 5, below.

The derivative rule also immediately also implies a novel approach to Stieltjes inversion as, at the quadrature points x[k],

$$\rho(x[k]) = w[k]/x'[k] \qquad (15b)$$

which is the subject of the Sections 3 and 4, as it suggests a straightforward manner of construction of ρ(x[k]) directly from the quadrature weights and points. The possibility of construction of the weight function from this quadrature information was well known to Markov, Stieltjes, and Chebyshev, and examples of their approach appear at the beginning of Section 4. Perhaps surprisingly the derivative rule leads to a simpler and far more efficient numerical method than these earlier approaches, an idea going well beyond that originally suggested by Heller, Ref. (21). At this point the result of Eqn. (15a,b) is referred to as the *derivative rule conjecture*.

### C) Special Cases: Derivative Rule May Sometimes Be Implemented Exactly

Taking advantage of the properties of the Chebyshev polynomials summarized the Appendix, a single example, also appearing in Refs. (18, 22), is illustrated. Consider the quadrature approximation for the integral

$$\int_{-1}^{1} g(x) \sqrt{(1-x^2)}\, dx. \tag{16}$$

The appropriate OPs are the Chebyshev polynomials of the second kind on the compact interval [-1,1]. These have, for a quadrature of order n:

$$\rho(x) = \sqrt{(1-x^2)} \tag{17a}$$

$$x[k] = x[k,n] = = -\cos[\, k\, \pi/(n+1)], \qquad k = 1,2,\ldots n \tag{17b}$$

$$w[k] = w[k,n] = \pi/(n+1)\, \sin^2[\, k\, \pi/(n+1)], \quad k = 1,2,\ldots n \tag{17c}$$

at the zeros of the $n^{th}$ degree polynomial the weight function $\rho(x)$ becomes

$$\rho(x[k]) = \sin[k\, \pi/(n+1)], \tag{17d}$$

where, as befits a non-negative $\rho(x)$ we have chosen the positive root, and inserted a "-" sign in Eqn. (16b) so that the abscissas increase as a function of "k" on the interval [-1,+1]. This would be immaterial for a Gauss quadrature approximation, but as we plan to take the derivative dx[k]/dk a non-negative result is chosen.

Consider now

$$x'[k] = dx[i]/di\,|_{i=k} = \pi/(n+1)\, \sin[\, k\, \pi/(n+1)], \tag{17e}$$

which may be compared to the ratio $w[i]/\rho(x[i])$ which is, from Eqns. (17c) and (17d), which also gives exactly x'[k].

The derivative rule conjecture is thus validated for the Chebyshev OP system of the second kind. Similar derivations, indicated in the Appendix, show validity of the derivative rule for all four sets of Chebyshev OP systems, as may also be validated for variations derived from them via linear transformations of the interval [-1,1], and or truncation thereof.

It is also useful to note, for later use, that

$$\lim_{n\to\infty} x'[k,n]) = (\pi/n)\sqrt{(1-x[k,n]^2)}, \qquad (18)$$

the later expression appearing, often, in the literature relating to universality and random matrix theory, for OP systems defined below. See Refs. (23-27), and references therein. The Eqn. (A.5) of the Appendix show, that both x[k,n] and x'[k,n] are indeed identical (i.e. universal) functions for all four sets of Chebyshev polynomials as n→∞, except for k near the boundaries, 1 or n.

### D) Initial Conclusions

This analytically exact results, for x'[k,n], above, and those in the Appendix, show that when x'[k] may be computed exactly the derivative rule conjecture is confirmed.

The conjecture may be morphed towards becoming a theorem when embedded in the concepts of universality of eigenvalue distributions for classes of random matrices. A crucial step (see Askey and Ismail, Ref. (7)) was the realization that the use of three term relations for the *orthonormal polynomials*, rather that of Eqn. (3) for the traditionally normalized OPs, provided a necessary systematization for considering families of related OPs and their perturbations, as well as being the formulation most closely connected to modern numerical calculations of the quadrature abscissas and weights, as in Refs. (Temme, DLMF, ref 1b)), Gautschi (3), Golub and Meurant (3), Gordon (5), all of which are based on a real symmetric tri-diagonal Jacobi or J-matrix.

It is the purpose of the next Section to indicate the outline of a proof of the conjecture, using techniques, now standard, but which did not exist in the 1970s.

## 3. Derivation of the Derivative Rule in the Context of Universality

### A) The J-matrix and Overview of the Derivation

The recursion relation for orthonormal (rather than simply orthogonal as for those defined in Eqn. (3)) polynomials, $p_n(x)$, is conveniently expressed in terms of the J-matrix coefficients $a_i$, $b_i$ (& noting that there are *many* conventions for labelling these coefficients), we follow the notation of Levin and Lubinsky, Ref. (23)

$$x\, p_n = a_{n+1}\, p_{n+1}(x) + b_n\, p_n(x) + a_n\, p_{n-1}(x), \quad n = 1,2\ldots \qquad (19)$$

initialized as $p_0(x) = 1/\sqrt{\mu_0}$, $p_1(x) = (x - b_0)/(a_1\sqrt{\mu_0})$. Gautschi, Ref. (3), p 29, lists these recursion coefficients for the classical OP systems in an alternative notation.

The polynomials to be considered to establish the validity of the derivative rule are of the Nevai-Blumenthal (N-B, hereafter) class, Nevai *et. al.*, Ref. (24), defined by

$$\lim_{n\to\infty} a_n = \frac{1}{2} \quad \text{and} \quad \lim_{n\to\infty} b_n = 0, \tag{20}$$

which includes, among others, the Chebyshev, Legendre, Gegenbauer, and the Repulsive and Attractive Coulomb Pollaczek polynomials; but, not classical OPs such as the Hermite and Laguerre systems. The first four of these are, in addition, special N-B systems having zeros in the interval [-1,1] and only an absolutely continuous (*a.c.*) Stieltjes measure on that interval. The Attractive Coulomb Pollaczek polynomials (defined below) have discrete mass points (i.e. isolated jumps in μ(x)), and zeros, for x < -1, and accumulating at x = -1⁻, in addition to the *a.c.* measure on [-1,1], as follows from their relationship to a Coulomb Schrödinger operator as in Ref. (18). See also Refs. (24), where it will be seen that this possibility was clearly foreseen by Blumenthal.

The J-matrix itself, with the conventions of Eqn. (19), is

$$\mathbf{J_n} = \begin{pmatrix} b_1 & a_1 & & & & 0 \\ a_1 & b_2 & a_2 & & & \\ & a_2 & b_3 & a_3 & & \\ & \ddots & \ddots & \ddots & & \\ & & & a_{n-2} & b_{n-1} & a_{n-1} \\ 0 & & & & a_{n-1} & b_n \end{pmatrix}$$

$$\tag{21}$$

and is an important and useful quantity, as its eigenvalues are the zeros of the $n^{th}$ order OP generated by the recursion coefficients, and the quadrature weights may be also easily constructed via the J-matrix, see Temme, Ref. (1b), Golub and Meurant, Ref. (3), for alternate approaches.

These N-B class polynomials have *universal* distributions of zeros on [-1,1] which follow a *clock-rule limit*. See Refs. (24-27). The term *clock-rule* first arose in discussions of the zeros of OP systems orthogonal on the unit circle, and thus the zeros might be thought of as numerals surrounding a clock-face. Simon, Ref. (26) section 8.4, provides a most informative set of illustrations of said clock-faces.

Using this universality, the derivation of the derivative rule consists of the following steps:

i) It is shown that the derivatives x'[k,n], x ∈[-1,1], follow the same universality as the clock-rule distribution of zeros. Said another way, this amounts to interpreting the spacings, x[k,n] - x[k-1,n] , of the clock-rule zeros of an order *n* polynomial, as approximations to those derivatives. As $n \to \infty$ this approximation becomes exact.

ii) It is then noted, and illustrated, that the equivalent quadrature ratios w[i]/ρ(x[i]), for the above N-B polynomials and for x ∈[-1,1], have the same universal behavior as x[k,n] - x[k-1,n]: namely on [-1,1], with the above assumptions the ratio w[k,n]/ρ(x[k,n]) is independent of the polynomial under consideration. This result, although stated differently and not so interpreted, has been derived by Levin and Lubinsky, Ref (23).

iii) Equating the two universal quantities x'[k] and w[k]/ρ(x[k]) provides a derivation of the derivative rule for the *a.c.* part of the measure.

In Section 4 the derivative rule will be shown to give an exponentially convergent solution of the Stieltjes inversion problem, not only for the above polynomials, taken as simple illustrative cases, but also for the Gegenbauer and Hermite OP systems. The latter of these is not in the N-B class, or the *regular* class of measures actually assumed in Ref. (23), indicating that the assumptions made in the derivation, as outlined above, are not necessary for the validity of the derivative rule.

**B) The Clock-Rule and Universality of the Derivatives x'[i].**

With the above assumptions of the N-B class, and assuming the ordering x[k,n]>x[k-1,n] for x ∈ [-1,1] the clock rule is stated, in Ref. (23), page 71, as

$$\lim_{n \to \infty} (x[k,n] - x[k-1,n]) \frac{n/\pi}{\sqrt{(1-x[k,n]^2)}} = 1. \qquad (22)$$

Using Eqn. (18 or A.5 from the Appendix, as the Chebyshevs are in the H-B class), above, this may be re-written as

$$\lim_{n \to \infty} \frac{(x[k,n] - x[k-1,n])}{x'[k,n]} = 1, \qquad (23)$$

Where x'[k,n] is the derivative with respect to k of the zero x[k,n] of any one of the Chebyshev polynomials. The fact that the ratio

$$\frac{(x[k,n] - x[k-1,n])}{x'[k,n]} = \frac{\Delta[k,n]}{x'[k,n]} \qquad (24)$$

tends to unity as n→∞, implies that the difference Δ[k,n] = x[k,n] - x[k-1,n] is indeed an approximation to that derivative. This convergence to the derivative is illustrated in Figures 1a,b,c. Where Δ[k,n] and x'[k,n] are plotted as a function of k, *for the first three types of Chebychev polynomials, for n = 10, 50, and 1000. The three Δ[k,n], (using x1[k,n], x2[k,n]…of equations of the Appendix ) are in black, and the three x'[k,n], (i.e. x1', x2', etc.) in red, and convergence of all six of these quantities to a common and universal curve is seen at *n* = 1000. The ratios of any two of these quantities, as a function of k, for fixed n would clearly be approaching unity, which is the clock-rule result.

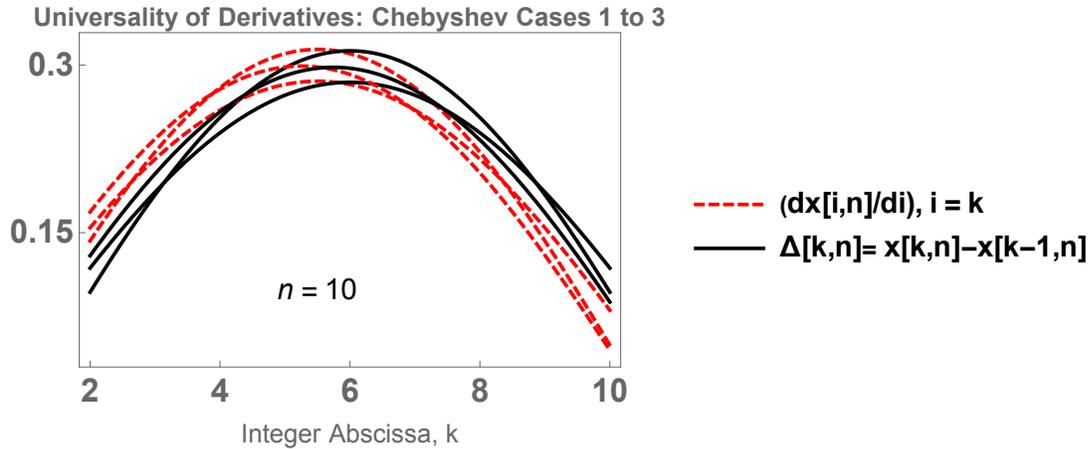

**Figure 1a.**

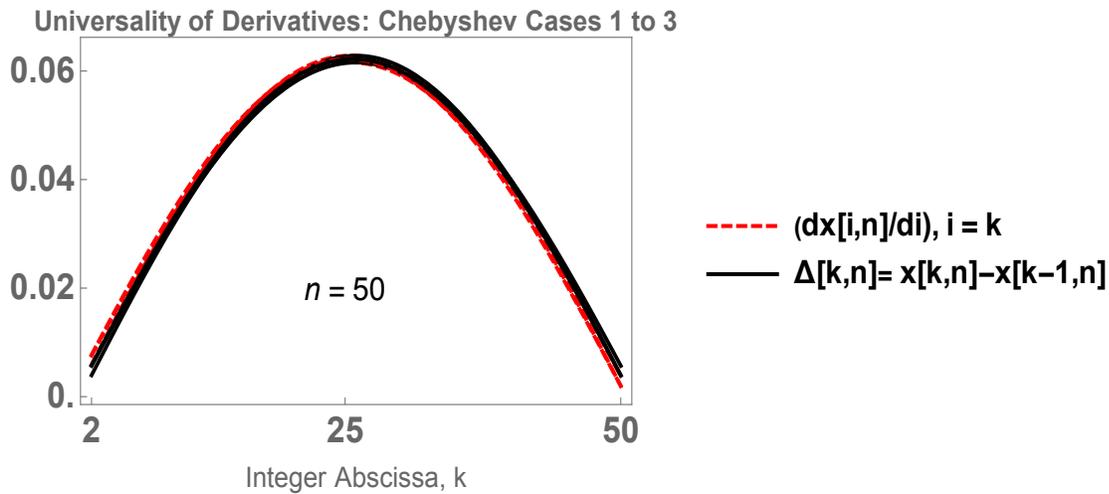

**Figure 1b.**

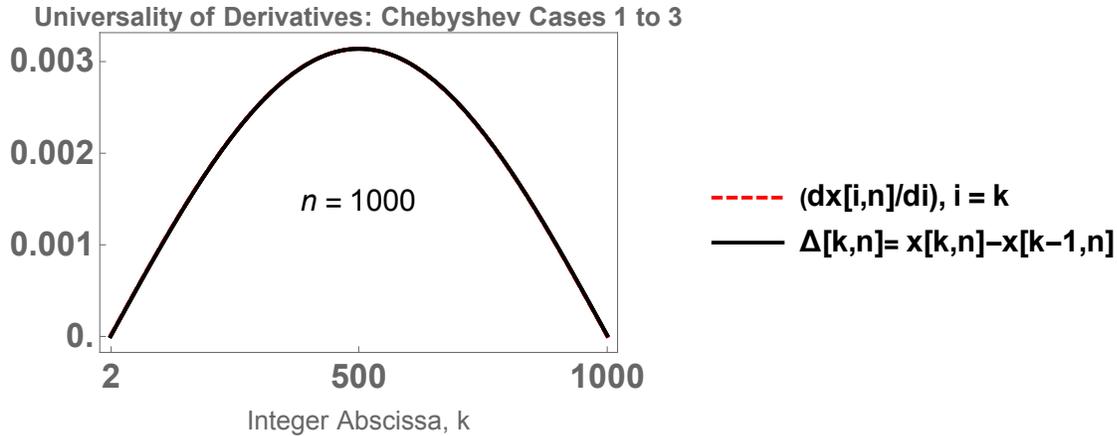

**Figure 1c.**
**Figure 1 Caption.** *The approach to universality is illustrated for the first three types of Chebyshev polynomials. See the Appendix. The data for the individual polynomials are not separately indicated as all converge identically. What is shown, for each polynomial, are the clock rule offsets of consecutive zeros $\Delta x[k,n] = x[k,n]-x[k-1,n]$, and the derivative rule, $x'[k,n]$ as evaluated at the first of these. The six graphs are converged, to visual accuracy, at n =1000, and thus their rations are approaching the common value of "one." The clock rule spacings are thus simply the derivatives $dx[k]/dk$ at the integer values of k as $n \to \infty$.*

Such is universality, as is clearly shown in Figure 1. This same factor,
$(\pi/n)\sqrt{(1-x[k,n]^2)}$ , arising here as a simple derivative, appears as the *weak-limit* of a canonical probability distribution for polynomial zeros in the usual, and more mathematically oriented, discussions of the clock-rule and its relationship to the zeros of random polynomials of the N-B class, see Simon, Ref. (25) page 215, and Ref. (27).

  Figure 1 which was designed to provide a clear illustration of the limit of Eqn. (24) does so, albeit with slowly convergent numerical results. The major source of error, and slow convergence, being the fact that the that the *backward differences* x[k,n] - x[k-1,n] while, indeed, giving an approximation to the derivative x'[k,n], are actually a rather poorly converged numerical estimate. This is well illustrated in Figure 2 where x'[k,n] is replaced by x'[k-1/2,n], with respect to which x[k,n] - x[k-1,n] provides a *central difference* approximation to the derivative, and thus showing a rather better numerical convergence, even for *n*=10. See Ref. (28) for a discussion of these different numerical approximations. The remaining differences are now actual differences in the derivatives of the different Chebyshev polynomial zeros, as seen for small n. As seen in Figs. (1a,b,c) these differences rapidly disappear.

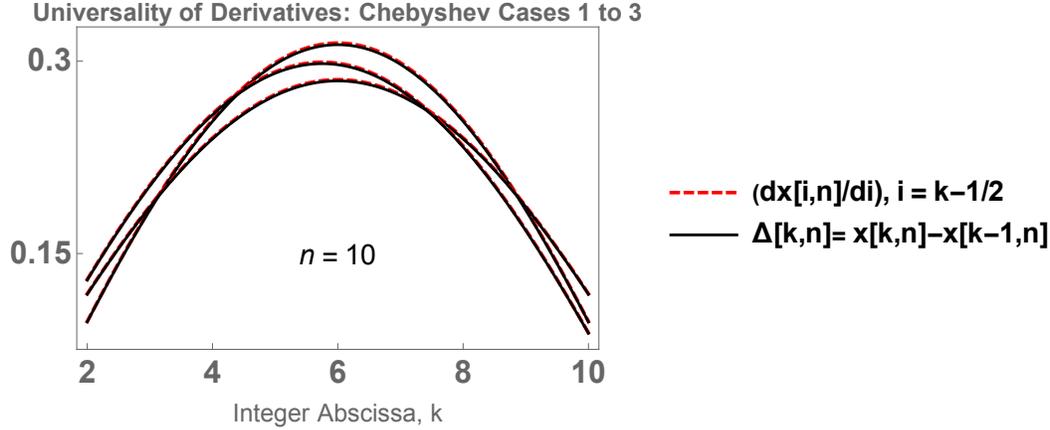

**Figure 2.** *A central difference approximation to the Δx[k,n] = x[k,n]-x[k-1,n], namely x'[k=1/2,n] gives more rapid convergence, that that shown in Figure 1. As n → ∞ the results are unchanged.*

A slightly more case arises for the Repulsive and Attractive Coulomb Polleczek polynomials, first considered by Yamani and Reinhardt, Ref. (18), and elaborated on in Refs. (29). Here, as these are less familiar, the Jacobi matrix elements (or J-matrix recursion coefficients) are:

$$a_n = \frac{1}{2}\sqrt{\frac{n(n+2l+1)}{(n+l+1+\frac{2Z}{\lambda})(n+l+\frac{2Z}{\lambda})}} \qquad (25a)$$

$$b_n = \frac{2Z/\lambda}{n+l+2Z/\lambda} \qquad (25b)$$

Here n = 1,2,…, as per the recursion of Eqn. (19), *l* is the angular momentum quantum number, *l* = 0,1,2…, and Z the product of the signed charges on the two interacting particles (in atomic units, where Z = ±1, ±2…), and λ is a scale parameter used to specify an $L^2$ Laguerre type basis for *discretization* of the Coulomb Schrödinger operator for the quantum Coulomb problem, see Ref. (18). λ must be chosen so that the Favard's Theorem inequality (*n* + *l* + 2Z/λ) > 0 is satisfied, which it may well not be for Z < 0, or the relation between the discretized Schrödinger operator and the Pollaczek polynomials breaks down, see also Ref. (18). While the recursion coefficients of Eqn. (25) are of the N-B class, and the *a.c.* spectrum exists only for x ∈ [-1,1] for both positive and negative Z.

The attractive case has, in addition, a discrete point spectrum for x < -1, corresponding to the bound states of the atomic system under consideration, and following the initial insight of Blumenthal, as confirmed by Nevai, Ref.(24), has an accumulation point at -1. Although, in this case, that is more quickly deduced as an actuality, rather than simply a possibility consistent with Eqn. (20), from the well-known spectrum of the Coulomb Schrödinger operator, Ref.(18,29). Thus, as labelling the polynomial zeros by an integer fails, for Z< 0, i.e. the Attractive case, to distinguish between the zeros in the discrete and absolutely continuous parts of the spectrum, the graphic displays of universal behavior are modified to an abscissa labelled by x ∈ [-1,1], rather than an integer, k. The energy differences are, however, indicated as being between consecutive zeros, both in the *a.c.* spectrum, as x[k,n] - x[k-1,n], without specification of the value of k.

Figure 3 illustrates the universality of the consecutive zeros of Chebyshev polynomials of types 1,2, and 3, and for both the Attractive (negative Z) and Repulsive (positive Z) Coulomb Pollaczek polynomials for Z = ±1, $l$ = 0, λ=4, and *n* = 2000. The Chebyshev zeros are known analytically, as before, but for the first time herein, the zeros of the polynomials under consideration have been calculated numerically by direct diagonalization of the J-matrix with working precision of 160 decimal digits. For the attractive, Z = -1, λ=4, case the first 40, out of 2000, zeros are outside the range of the *a.c.* spectrum, and are not shown. Only ~1% of the data is shown, and a shift added, to the values of x shown for the AC and RC Pollaczek polynomials, to allow clear distinction to be made. This same sampling is used in Figure 4.

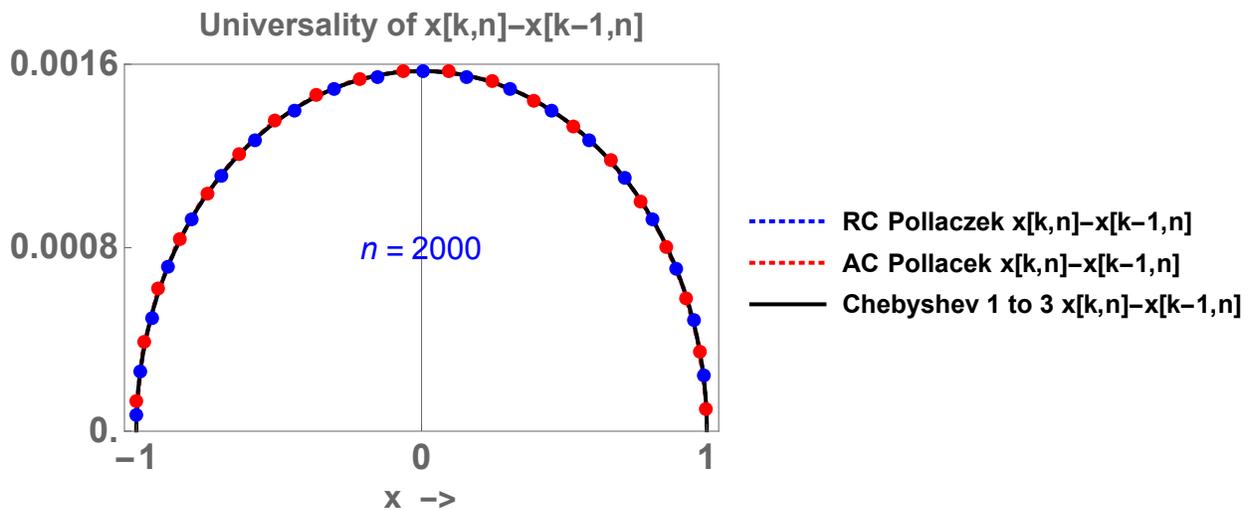

**Figure 3.** *The abscissa differences Δx[k,n] = x[k,n]-x[k-1,n], for the first three Chebyshev polynomials (solid line with all values of k displayed) are compared to those (with x ∈ [a,b]) for the Attractive and Repulsive Coulomb Pollaczek polynomials. The data for these latter two being only shown for every 100^th value of k, with an offset of 50, so that they may be easily seen and differentiated. All five sets display the same universal behavior, not surprising as all are in the N-B class of OP systems.*

## C) The Derivative Rule, and Universality of w[k]/ρ(x[k])

Following Levin and Lubinsky, Ref. (23), p 78, and the above discussion, we have assuming OPs in the N-B class, in the limit of large *n*, and assuming that any "jumps" in µ(x) are outside of [-1,1],

$$\frac{(x[k,n]-x[k-1,n])}{x'[k,n]} = (x[k,n] - x[k-1,n])\widehat{K}_n(x[k,n], x[k,n]). \quad (26)$$

Where $\widehat{K}_n(x[k,n], x[k,n]) = K_n(x[k,n], x[k,n])\rho(x[k,n])$, $K_n(x[k,n], x[k,n])$ being the usual reproducing kernel (see Ismail, Ref. (2) sections 2.2 and 2.4) which takes the value 1/w[k,n] for the arguments as shown. Comparison of the two sides of Eqn. (26) gives, immediately,

$$x'[k,n] = w[k,n]/\rho(x[k,n]) \quad (27)$$

which is the derivative rule. Further, the quantity w[k,n]/ρ(x[k,n]) is seen to have the same universality as $(x[k,n] - x[k-1,n])$, or x'[k,n], for large n.

Namely for OPs in the N-B class w[k,n]/ρ(x[k,n]) is a universal function for an *a.c.* spectrum on [-1,1] as n→ ∞. This is illustrated in Figure 4 for the Chebyshevs, Legendre, and the Attractive and Repulsive Coulomb Pollaczek polynomials, these latter three being calculated numerically as in Figure 3, with the w[k,n] calculated using the algorithm of Golub and Meurant, Ref (3), with a working precision of 160, again for n = 2000. The Pollaczek parameters, and sampling, for the Pollaczek and the Legendre, polynomials are the same as for Fig. 3, and listed in the paragraph describing that Figure.

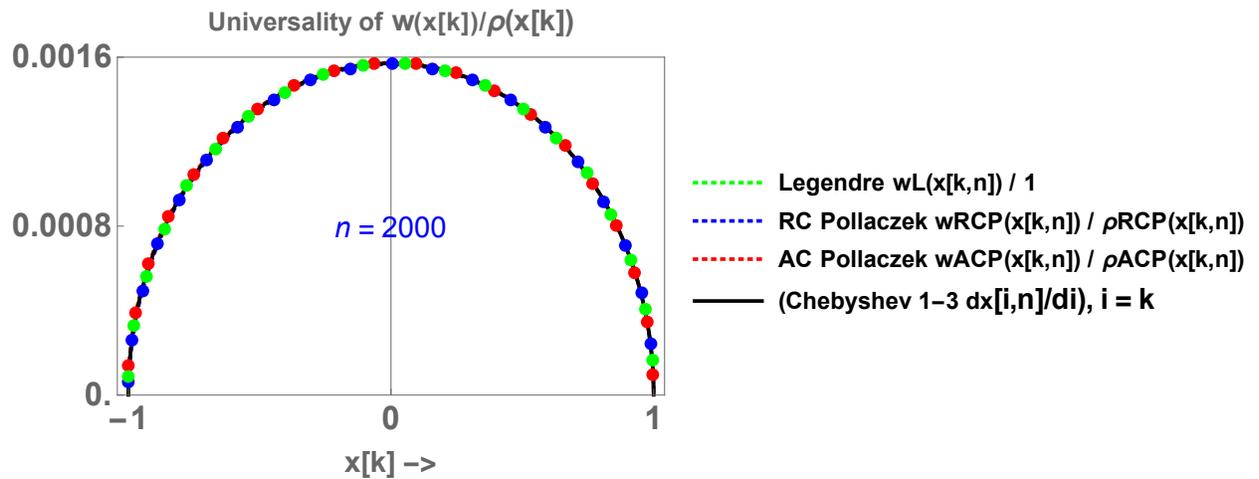

**Figure 4.** *The values of x'[k,n] = w[k,n,]/ρ(x[k,n]) for the first three Chebyshev polynomials (solid line with all values of k displayed) are compared to those (with x ∈ [a,b]) for the*

*Legendre, and the Attractive and Repulsive Coulomb Pollaczek polynomials. The data for these latter three being only shown for every 100th value of k, with offsets of 30 and 60, so that they may be easily seen and differentiated. All s sets display the same universal behavior.*

The result is particularly striking given the differences in the ρ(x) weight functions themselves as illustrated in Fig. 5, noting that any differences in normalization of the weight functions cancels in the derivative rule ratio illustrated in Fig. 4.

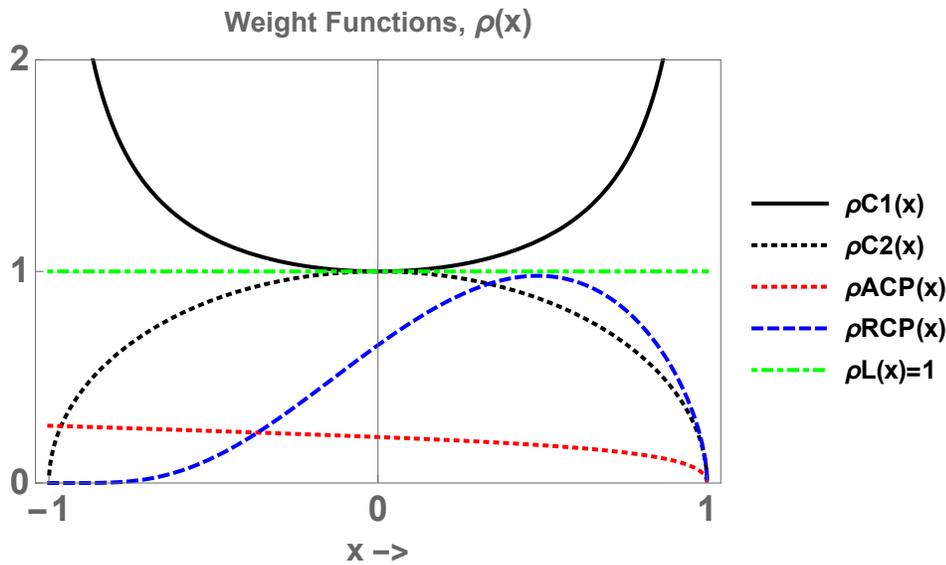

**Figure 5.** *The weight functions, ρ(x), for type 1 and 2 Chebyshev, the Legendre, and the Attractive (AC) and Repulsive (RC) Coulomb Pollaczek polynomials. These are not universal functions, although all of these OP systems are of the N-B class. As shown in Figure 4, the ratios w[k,n]/ρ(x[k,n]) are universal functions of x[k,n] ∈ [-1,1]. This latter restriction on x being important as the weight function for the Attractive Coulomb Pollaczek system has positive discrete (positive) jumps for x < -1, which accumulate at -1⁻. These are not shown, being outside [-1,1]. This leads to quite different areas of the a.c. parts of the Attractive and Repulsive weight functions.*

The formulae for the Chebyshev weight functions are in the Appendix, the Legendre weight function is ρ(x)=1, and the Coulomb Polleczek (CP) weight function, x ∈ [-1,1], is, see Refs. (18, 29),

$$\rho^{CP}(x) = \frac{2^{(2l+1)}}{\pi} e^{-(2\theta(x)-\pi)\gamma(x)} (1-x^2)^{(l+1-1/2)} |\Gamma(l+1+i\gamma(x))|^2 \frac{(l+1+2Z/\lambda)}{\Gamma(2l+2)}, \qquad (28a)$$

where

$$\theta(x) = \text{Arccos}(x), \quad \gamma(x) = Z/\kappa(x), \quad \text{and} \quad \kappa(x) = \sqrt{\frac{\lambda^2(1+x)}{4(1-x)}}, \quad (28b)$$

with essential singularities at each end of the *a.c.* interval. It is worth noting that unless Favard's inequality, $(l+1 + 2Z/\lambda) > 0$, is satisfied for Z <0, $\rho^{CP}(x)$ will be zero or negative. Thus, the importance of the inequality, which comes into play, as for example, in the quantum theory of the Hydrogen atom ground state for $l = 0$, and where Z = -1, and thus the requirement that $\lambda > 2$. Note that there is no fundamental physics which determines $\lambda$, it simply being a non-linear variational parameter in an OP-type basis, as discussed in Ref. (18).

The overall Coulomb Pollaczek weight functions, for Z positive or negative are both normalized to unity. The discrepancy, clearly visible in Fig. 5, between the areas under the Attractive and Repulsive weight functions as shown, is due to the additional discrete weights for Z < 1. The discrete weight function, corresponding to the quantum bound states, for x < -1, is discussed by Bank and Ismail, Ref. (29) and Ismail, Ref. (2), section 5.8. Where Z, z in their notation, needs to be replaced by Z/$\lambda$, with $\lambda$ chosen appropriate to both the scale of the physical problem at hand, and to satisfy the Favard inequality. These discrete weights are not part of the *a.c.* spectrum as shown on [-1,1], nor included in the above weight function. The numerically determined quadrature weights, w[k], corresponding to the zeros x[k] < -1 sum to ~0.58787, consistent with the "missing" area under the *a.c.* weight function.

## 4. Applications of the Derivative Rule to Stieltjes Inversion

The question at hand is: given the w[i] and x[i], how might we best approximate the weight function $\rho(x)$ for x $\in$ [a,b]? This is explored here by consideration of the model problem of determination of the weight functions, $\rho(x)$, for the typical members of the classical OPs.

Stieltjes inversion by construction of a step-histogram approximation to $\mu(x)$, based on consecutive partial sums of the w[i], and via the Markov-Chebyshev inequalities guaranteeing convergence, goes back into the 19th century as discussed in Ref. (4). Chihara, Ref. (2) has provided a clear outline of this methodology, with references to the historical material, and his approach is followed here. This traditional method will be shown, in Sub-Section A to give results for $\rho(x)$ which converge as ~$1/N^2$, i.e. to modest precisions for computationally easily obtainable values of *n*. However, to obtain results of high precision rather large values of n are necessary: for example, convergence of $\rho(x)$ to 6 or 7 decimal digits requires N ~ 2000. For many purposes the histogram method is sufficient, especially when the input is experimental data, and these technique has been exploited in atomic physics by Langhoff, and Langhoff et. al., Ref. (30); and by Gordon, Ref. (5) who obtains upper and lower bounds of improved precision to quantities of interest in statistical physics. All three of these contain extensive referencing of the original literatures.

In Sub-section B) it is seen that, rather than power-law convergence, use of the derivative rule gives exponential convergence of ρ(x) for the classical and N-B polynomials. Thus, for example $n = 20$ might give precisions of a part in $10^{20}$, and $n = 40$ a part in $10^{40}$, thus opening the way to a completely new set of possibilities. This is illustrated for Chebyshev and Gegenbauer systems.

In Sub-section C) the case of the Hermite OPs is briefly discussed. These are classical, but not N-B, systems, and, further, have weight functions on an open, and infinite interval. It is seen that exponential convergence on the open interval $(-\infty, \infty)$ is obtained into regions where ρ(x) is of order $10^{-200}$. Brief comments on the more difficult case of determining the Laguerre weight function are included.

### A) Traditional Histogram Methods for Stieltjes Inversion

Beginning with tradition, consider the OP system by OPs orthogonal on [a,b] with weight function dμ(x): Define $\mu_n(x)$ as a bounded, right continuous, non-decreasing step function whose steps are at the x[i,n], and of positive magnitude w[i,n]. Using the Heaviside relation between such step functions and the Dirac-delta function, δ(x-x'), (see, for example Ref. (5)) this is most simply expressed as:

$$d\mu_n(x) = \sum_{i=1}^{n} w[i,n]\delta(x - x[i,n]) \qquad (29)$$

This is equivalent to

$$\mu_n(x) = \sum_{i=1}^{n} w[i,n]\,\Theta(x - x[i,n]) = \int_a^x d\mu_n(x) \qquad (30)$$

Θ(x –x[i]) being the unit Heaviside step function, 0 for x < x[i], 1 for x ≥ x[i]. Thus, omiting the n's, for x[i] ≤ x< x[i+1], i+1 ≤n:

μₙ(x) = 0                                         a ≤ x < x[1],

$\mu_n(x) = \sum_{j=1}^{i} w[j]$                 x[i] ≤ x <x[i+1],     1 ≤ i <n,       (31)

$\mu_n(x) = \sum_{j=1}^{n} w[j] = \int_a^b d\mu(x)$     x[n] ≤ x <b,

which gives a histogram. This is the Chihara *ansatz* of Ref.(2), page 56, which may be traced back to Markov, Chebyshev, and Stieltjes, as in the discussion and references in Section 2Biii. For x in the absolutely continuous spectrum,

$$\rho(x) = \frac{d}{dx} \int_a^x \mu_n(x) dx \qquad (32)$$

and Chihara ansatz of Eqns. (31) can be smoothly interpolated, between the x[i], to allow implementation of the derivative in Eqn. (32). Carried out directly, however, this leads only to convergence as 1/N.

    A better numerical approach is to smoothly interpolate $\mu_n(x)$ at the midpoints (x[k]+x[k+1])/2, and then take the derivative of Eqn. (32) of this interpolation, as in Ref. (31), at the points x[k]. This interpolation, which Langhoff calls the Stieltjes derivative, see Ref. (30), gives approximations to $\rho(x)$ which converge as $1/N^2$, as now illustrated. Taking as a simple example the Chebyshev polynomials of the second kind, differentiation of the numerically interpolated $\mu_n(x)$ gives the results of Figures 6a, and b. The error, as shown, is simply $\rho^{Exact}(x) - \rho^{Approx}(x)$ over the range of the *a.c.* spectrum, $x \in [-1,1]$. The approximate error indeed scales as $\sim 1/N^2$, over the range N = 100 to 200,000. The fact that x[i] and w[i], are known exactly, for arbitrary N, for these polynomials has been exploited, as for N= 200,000, an unlikely set of data had one to diagonalize a J-matrix of dimensions $N \times N$. Yet the error, near x[k] = 0, is still $\sim 10^{-10}$, and completely insensitive to the interpolation order. Thus, no computations of this type have been attempted, for other members of the N-B class of OPs where the quadrature abscissas and weights are not known analytically. Ref. (30) contains many fully worked examples, where a far lower precision is actually necessary or sought.

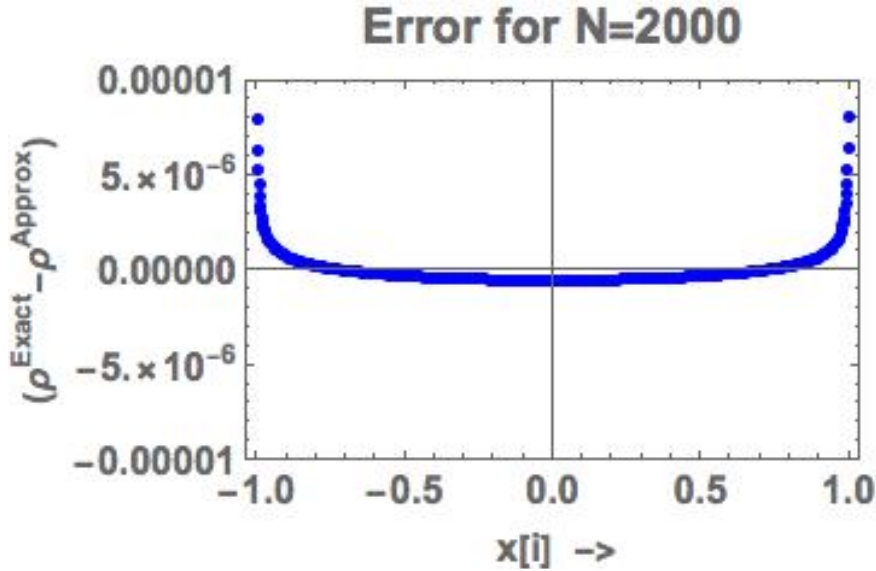

**Figure 6a.**

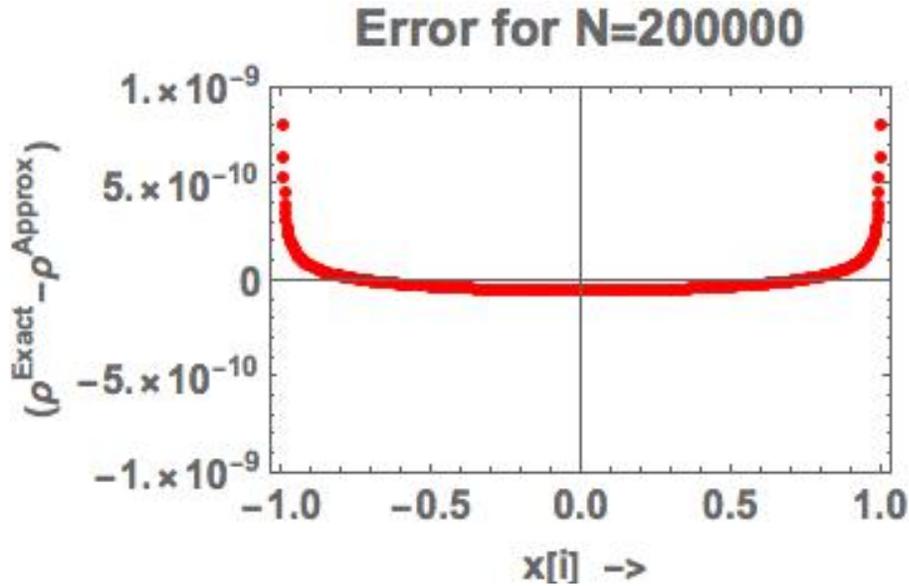

**Figure 6b.**

**Captions for Figures 6a,b:** *Differences between exact and approximate Chebyshev weight functions of the second kind, the approximation being the 'symmetric' version of the Chihara approximations, appropriately numerically interpolated with order ~10, see Ref. (31). In 6a, N = 2000; for 6b, N= 200000. The convergence is (approximately) of a power law type, convergence as ~1/N², and thus rather too slowly for convenient numerical work using standard methods of obtaining the x[i] and w[i] by diagonalization of an N×N J-matrix. Note that although analytic derivatives of x[i] are available, no use of that fact has been made in the above results.*

## B) The Derivative Rule for Stieltjes Inversion: Chebyshev Polynomials

Figures 7a,b,c give examples for N = 20, 40, and 200 of the efficacy of the derivative rule, Eqn. (27). Here the only input is the known fixed n weights, w[k] and quadrature abscissas, x[k] for integral values of $k$. The derivatives are obtained via numerical interpolation of the x[k], as a function of $k$, with direct differentiation of the resulting expression, and then evaluation of the derivatives at the integer points $k$. With use of the derivative rule convergence is to one part in $\sim 10^{170}$ for N = 200. This may be compared to the results of the previous Sub-section, where for N = 200000, convergence was one part in $\sim 10^{10}$

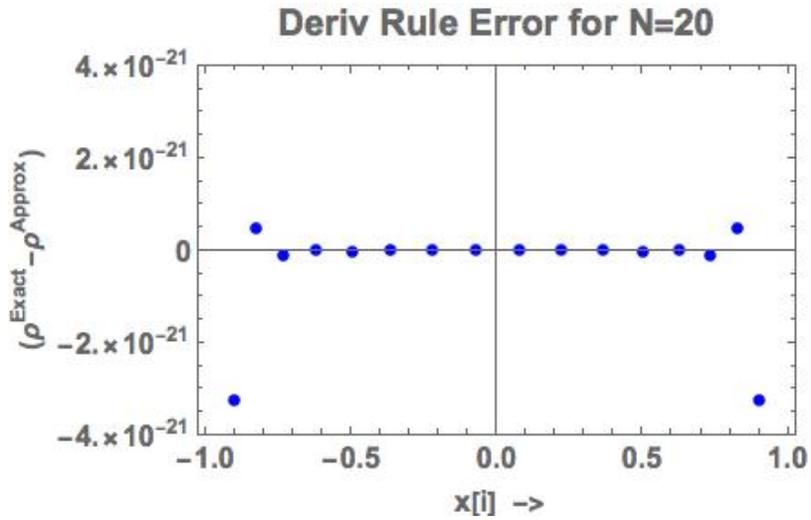

**Figure 7a.**

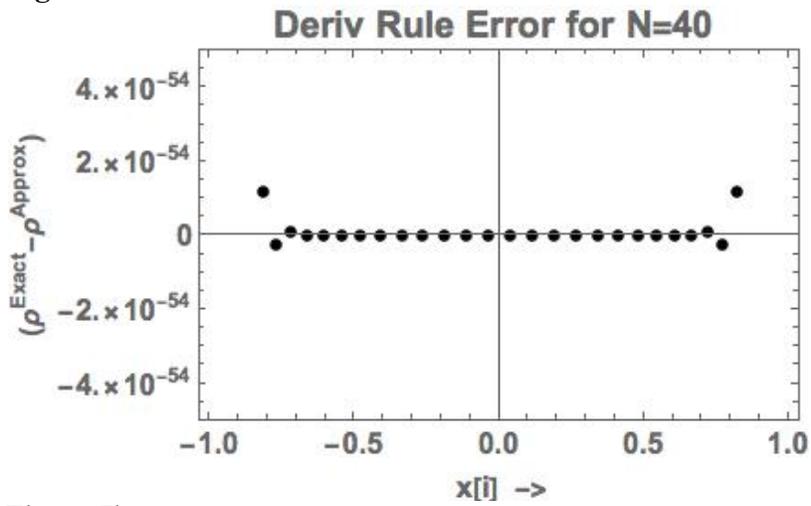

Figure 7b

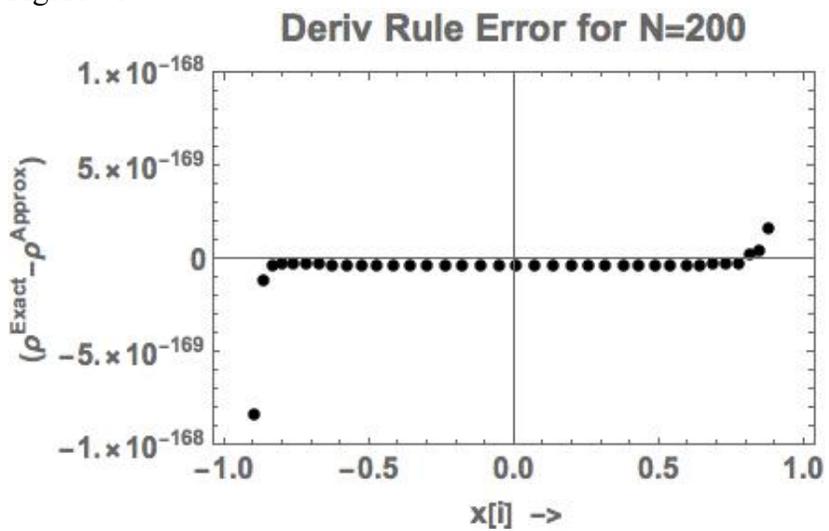

Figure 7c.

**Caption for Figure 7a,b,c.** *Here we show "typical, i.e. not optimized" convergence of derivative rule estimates of the Chebyshev weight function to $\rho(x) = (1-x^2)^{1/2}$, expressed as simply the difference between the approximate and exact functions. As N increases from 20 to 200, and an interpolation order of ~10, we see (at least for the central parts of the weight function) an "improvement" of nearly 150 orders of magnitude. In fact looking at the three graphs, a,b,c, and in Table 1, the convergence is of the form $10^{-\alpha N}$, with $\alpha$ the order of 1, signifying at least approximate exponential convergence, rather than the algebraic, $1/N^2$, convergence of the corresponding Figures 6a,b. Note that for Figures 7a,b 100 digit precision arithmetic was used, while, of necessity for 3c, 200 digits were carried throughout the calculation, see Ref (32).*

Table 1, below, shows that one part in $10^{170}$ is not an accident: the convergence is better than exponential in n. Why is this of note? As Charles Schwartz, Ref. (33), might have said as we are paraphrasing his analysis of the rate of convergence of fitting trial functions to desired functions in variational calculations, noting that exponential convergence is a special case: "We, at last, have a series of approximations *of the right shape*. Namely, boundary conditions, and points of non-analyticity, if any, have been properly included in our approximation scheme. Any of these, if mismatched, will produce only power law convergence, and with the stronger the mismatch the lower the degree of that power law convergence."

**Table 1.** *Convergence illustrating the better than exponential convergence of derivative rule computations of the Chebyshev weight function $\rho(x) = (1-x^2)^{1/2}$, near its maximum, at x = 0, as a function of N, via the derivative rule, and via an Hermite interpolation order N-1, of the abscissas x[k,n] to obtain x´[k,n]. In this case the interpolation order\* was pre-set to its maximum value N-1, and the fact that the weights and abscissas are known exactly, rather than needing to be calculated numerically, was used to full advantage, leading to results better than to be expected in more typical calculations, where these would be calculated numerically, and thus known only to finite precision. The reported error is at x[N/2] = ~ 0, namely near the maximum of the exact and approximate weight functions, and thus near the center of the interpolation range, likely enhancing the results. Figures 7a,b,c are perhaps better examples, as no particular optimizations or choice of points to make detailed comparisons have been chosen.*

| Value of N | Interpolation Order* | Approximate Error (to nearest power of 10) |
|---|---|---|
| 10 | 9 | $10^{-10}$ |
| 15 | 14 | $10^{-15}$ |
| 20 | 19 | $10^{-25}$ |
| 40 | 39 | $10^{-58}$ |
| 60 | 59 | $10^{-99}$ |

\*The term *interpolation order*, is discussed in Ref. (31).

### C) The Derivative Rule: Gegenbauer Polynomials

As a second example we consider the $l = 20$ Gegenbauer (also known as ultrashperical) weight function $\rho(x) = (1-x^2)^{(l + 1/2)}$ on $x \in [-1,1]$. In this and case the $x[i]$ and $w[i]$ are now computed from numerical J-matrix diagonalizations of Temme, Ref. (1b) and Golub and Meurant, Ref. (3). These polynomials, as are the Chebyshevs, in the N-B class. Our goal is, again, to reconstruct the above $\rho(x)$ from this input. Figure 8 gives a low resolution comparison of the weight function construction for N = 500, a rather larger than needed value for most applications, but a value which allows in Figures 9, development of more useful investigations to how to assess errors from such computations, and then how to improve convergence near x at the boundaries, at $x = \pm 1$. Again the derivatives $x'[k]$ are calculated numerically by Hermite interpolation, Ref. (31), as a function of k, followed by analytic differentiation of the Hermite fit, evaluated at the integers i = 1,2,...n.

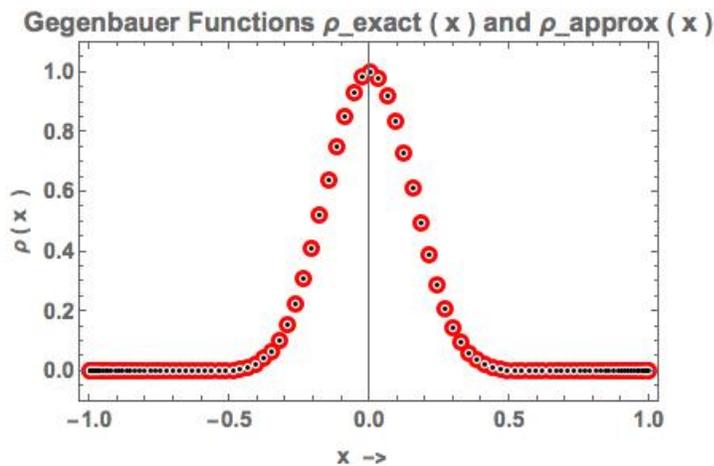

**Figure 8.** Exact, black dots, and approximate values, red circles, of the Gegenbauer l = 20 weight function for N = 500. Only ~10% of the data is displayed. The centers of the circles are at the approximated values of the weight function.

**Table 2.** *Approximate convergence of results for the l=20 Gegenbauer weight function, $\rho(x) = (1-x^2)^{(l + \frac{1}{2})}$ illustrating the approximately exponential convergence of derivative rule computations of $\rho(x)$, as a function of N (=n), calculated via the derivative rule, and via Hermite interpolations of the interpolation order* n = N-1, and not otherwise optimized, of the abscissas $x_N[k]$ to obtain $x'[k]$. Values of the errors shown are taken at the maximum of the weight function, i.e. at x = 0. In this case the exponential convergence is of the form $10^{(-N/\alpha)}$, and the values of $\alpha$ do differ, perhaps even systematically, although are approximately constant for $N \geq 61$. The broader nature of derivative rule approximation error is well illustrated in Fig. 9.*

Value of N     Approximate Error   Approx. Error as a Function of N
                                                                                 (showing power law convergence)

| | | |
|---|---|---|
| 11 | $10^{-5.1}$ | $10^{-N/2.15}$ |
| 21 | $10^{-8.7}$ | $10^{-N/2.43}$ |
| 41 | $10^{-15.4}$ | $10^{-N/2.66}$ |
| 61 | $10^{-21.8}$ | $10^{-N/2.79}$ |
| 101 | $10^{-34.4}$ | $10^{-N/2.94}$ |
| 201 | $10^{-65.3}$ | $10^{-N/3.08}$ |
| 401 | $10^{-126.3}$ | $10^{-N/3.17}$ |

*The term *interpolation order*, is discussed in Ref. (31).

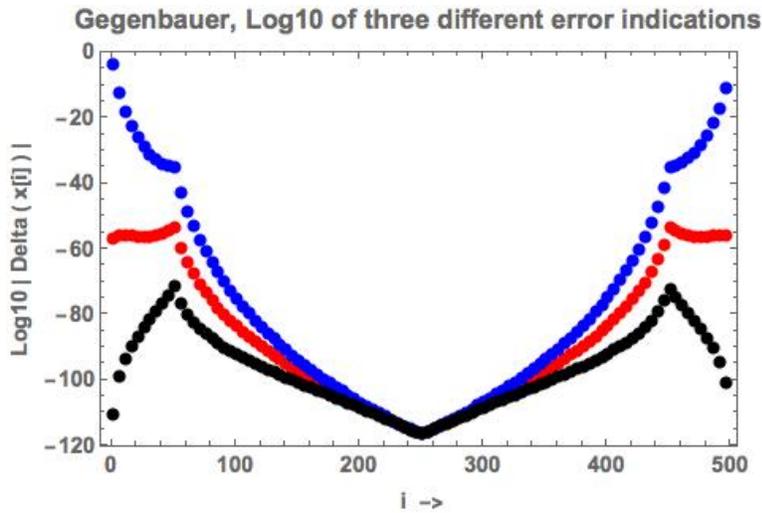

**Figure 9.** *As in Figure 7, except now differences between $\rho(x)\_exact$, and $\rho(x)\_approx.$, simply labeled as 'Delta,' as a function of $x[i]$, are displayed in three different ways, appropriate in different ranges of the value of $\rho(x)\_exact$. In this case N is fixed at 500, and the Hermite interpolation\* order fixed at 50, as we are illustrating different error presentations, rather than convergence as a function of N. The blue dots show $Log_{10}[\ |\rho\_exact - \rho\_approx| / \rho\_exact]$ at the abscissas labeled by k, running from 1 to 500. This measure of the error is most appropriate near the maximum of $\rho\_exact$. The red dots are simply $Log_{10}[\ |\rho\_exact - \rho\_approx|\ ]$ and appropriate for values of $\rho(x)\_exact$ in its intermediate, to large, ranges of value. For small $\rho(x)$ the black dots show the errors in the form $Log_{10}[\ |\rho\_exact - \rho\_approx| \times \rho\_exact]$, thus de-emphasizing the difference $|\rho\_exact - \rho\_approx|$ when $\rho\_exact$ is near zero.*

*The term *interpolation order*, is discussed in Ref. (31).

All of the measures of accuracy of Figure 9 clearly indicate poor convergence at the boundaries of the orthogonality interval [-1,1]. In the case here ρ_exact[x[1]] is the order of $10^{-57}$ out of a weight function normalized to 67282234305 /549755813888 = ~ 0.3845. This also clearly shows why [|ρ_exact – ρ_approx|/ ρ_exact] appears (via the blue dots) to be quite large at i = 1, even though, with magnitudes near $10^{-57}$, the exact and approximate values of the weight function are in agreement to about 4 (base 10) digits, which is actually, perhaps, worthy of some note.   Without going into detail, it should be added that use of point-wise continued fraction fits, which can then be analytically differentiated, see Refs. (10,12), rather than the Hermite interpolations of Mathematica, can easily add 8 to 12 additional digits near these boundaries should precisions at that level be warranted for quantities of this small magnitude. Yet another manner of displaying such data with very a broad range of values is used in Figures 10a,b where the weight function being approximated varies by as much as 14 to 200 orders of magnitude.

D) **Derivative Rule:  Inversion out of the Nevai-Blumenthal Class, Hermite Polynomials.**

The Laguerre and Hermite systems are important examples of non-N-B class OPs. The derivative rule, as originally stated, works well for both, again with exponential convergence away from boundaries. The Laguerres suffer from issues of instability of interpolation near the boundary x = 0 where the weight function is, in fact, at its maximum. These are largely handled by use of point-wise continued fraction fits to the x[i] near i = 1, as in Refs. (10,12), which may then be analytically differentiated. Perhaps this is a case where an appropriate alternative to the N-B universality could be of use.

With this latter not yet in hand, focus will be on the Hermite system which is characterized by orthogonality on the open interval (-∞,∞). Typical derivative rule inversions are shown in Figures 10a,b. for N = 21 and 401.  The boundary regions showing lower precision than at x = 0, but not inconsistent with the size of the weight function itself.

Shown are the Log10 of the weight function itself, $\rho(x) = e^{-x^2}$, at the abscissas x[i], and as compared to the Log10 of the difference between the actual value and its derivative rule estimate. For N = 21 the approximate value is accurate to ~8 digits at x[11]=0 where ρ(x) = 1, and to ~3 digits at x[1], (and x[21]), where ρ(x) = ~$10^{-14}$.  For n = 401, over the range x = ± 21,  ρ(x) has decreased to ~ $10^{-200}$ at the (open) boundaries, a smaller magnitude than would usually be of interest.   At x = ± 16,  ρ(x) ≅ $10^{-120}$, and the precision is ~70 digits.

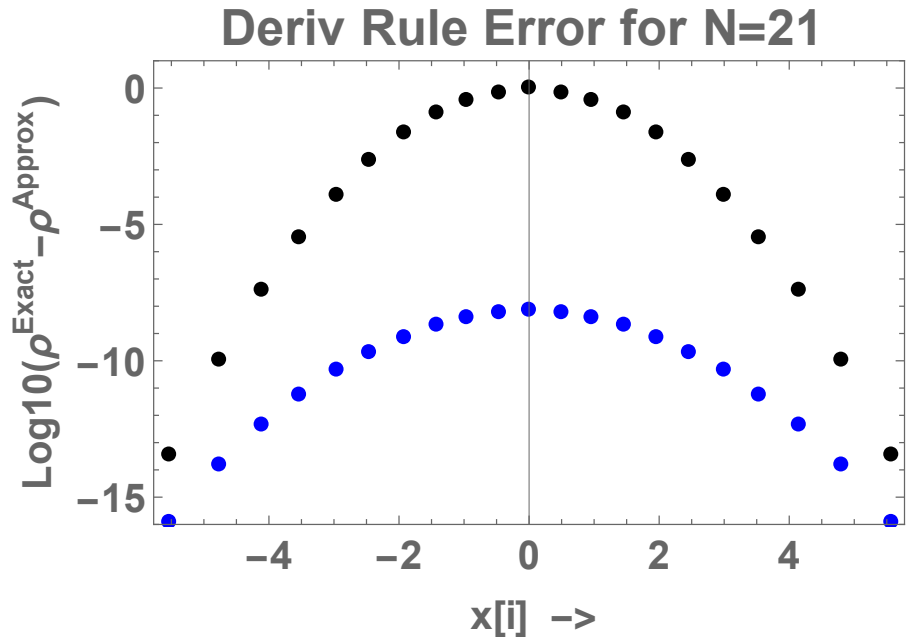

**Figure 10a.** *The upper black dots are Log10( ρ_exact(x)) to indicate the magnitude of the Hermite weight function, which varies over 14 orders of magnitude in this example. In comparison, the lower blue dots, indicate the error in the derivative rule computation, for N = 21, of the Hermite weight function as measured by Log10(|ρ_exact(x)) - ρ_approx(x)|).*

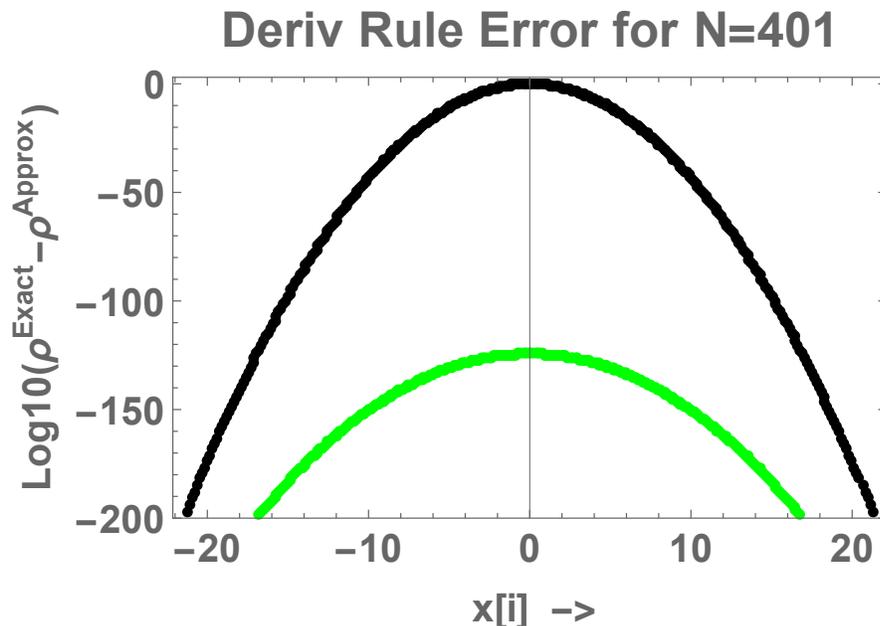

**Figure 10b.** *The upper black dots are Log10( ρ_exact(x)) to indicate the magnitude of the Hermite weight function itself, which varies over 200 orders of magnitude in this example. In comparison, the lower green dots, indicate in the error in the derivative rule computation of the Hermite weight function, for N = 401, as measured by Log10(|ρ_exact(x)) - ρ_approx(x)|). For N = 401 the dots have merged into what appear as solid lines. The calculations of the N =*

*401 results required a precision of 300 decimal digits, unusually high, but required by the nature of the weight function itself over this large range of x.*

Table 3 shows the power law convergence at x = 0, for the same values of N as for the N-B class Gegenbauer polynomials of Table 2. These results are quite intriguing as they are almost identical, suggesting that an underlying pattern has perhaps been exposed, yet to be understood. It is evident from this numerical work that the derivative rule is working well, outside of the N-B class where a derivation has now been supplied.

**Table 3.** *Approximate convergence of results for the Hermite weight function, Exp[-$x^2$], on the unbounded interval (-∞, ∞) illustrating the approximately exponential convergence of derivative rule computations of $\rho(x)$, as a function of N( = n), calculated via the derivative rule, and via Hermite interpolations of the interpolation order\* n = N-1, and not otherwise optimized, of the abscissas x[k,n] to obtain x′[k,n] . Values of the errors, simply $\rho$\_exact(x)-$\rho$\_approx(x), shown are taken at the maximum of the weight function, i.e. at x = 0. The exponential convergence is of the form $10^{\wedge}(-N/\alpha)$, and the values of $\alpha$ do differ, perhaps even systematically, although are approximately constant for N ≥ 61. Interestingly the errors as a function of N are almost identical to those of Table 2, for an l = 20 Gegenbauer weight function on the bounded interval [-1,1].*

| Value of N | Approximate Error (to nearest power of 10) | Approx. Error as a Function of N (showing power law convergence) |
|---|---|---|
| 11  | $10^{-5.5}$  | $10^{-N/2}$ |
| 21  | $10^{-8}$    | $10^{-N/2.6}$ |
| 41  | $10^{-15}$   | $10^{-N/2.8}$ |
| 61  | $10^{-21}$   | $10^{-N/2.9}$ |
| 101 | $10^{-34}$   | $10^{-N/3.0}$ |
| 201 | $10^{-64}$   | $10^{-N/3.1}$ |
| 401 | $10^{-127}$  | $10^{-N/3.2}$ |

\*The term *interpolation order*, is discussed in Ref. (31).

## 5. Origins of the Derivative Rule: Spectral Expansions of Schrödinger Operators

In the above analysis it has been assumed that knowledge of a three-term recursion relation implies the x[k] and w[k], the weight function then being constructed from the derivative rule. What if the three-term recursion relation is not at hand? What if the w[k] are not known? In the case of Schrödinger operators useful inversions are still possible, and this use actually gave rise to the first usage of what has been called the derivative rule in the above methodology for Stieltjes inversion. The course of this development may be overviewed in Refs. (12, 18, and 34 through 37) and references therein.

## A) "Direct Construction" of the Discontinuity in a Schrödinger Spectral Resolution from Two Discretized Approximations to the Resolvent

Suppose we wish to evaluate a quantum transition probability for a transition between an $L^2$ bound state $\varphi(x)$ to a continuum eigenstate, $\psi(E)$ of a Schrödinger Hamiltonian $\mathcal{H}$, where

$$\mathcal{H}\psi(E) = E\psi(E), \tag{33}$$

with $<\psi(E), \psi(E')> = \delta(E - E')$, and $\psi(E)$ being an appropriate non-$L^2$ (continuum) eigenfunction of $\mathcal{H}$, and $\delta$ the Dirac delta-function. Here $<,>$ is a scalar product in an appropriate Hilbert space, and to simplify the notation all functions are assumed real. Such a probability will be, assuming a proper choice of units, be proportional to

$$|<\varphi(x)t(x), \psi(E)>|^2, \tag{34}$$

with $t(x)$ being an appropriate transition operator. Often the non-$L^2$ eigenfunctions, $\psi(E)$, are difficult to work with, especially in a many-particle case, so we wish to find the values of matrix elements such as that of Eqn. (34) without the direct construction of the $\psi(E)$.

How is this question related to the Stieltjes inversion problem discussed in the previous Sections? Suppose that the *resolvent* $(z - \mathcal{H})^{-1}$ for complex z not in the spectrum of $\mathcal{H}$, has the spectral resolution (bound states being ignored):

$$(z - \mathcal{H})^{-1} = \int_0^\infty \psi(E) \frac{1}{z-E} \psi^*(E)\, dE, \tag{35}$$

then the matrix element $<\varphi t, (z - \mathcal{H})^{-1} \varphi t>$ would have the related spectral expansion:

$$\mathcal{F}(z) = \int_0^\infty \frac{|<\varphi(x)t(x),\psi(E)>|^2}{z-E} dE, \tag{36}$$

which takes things back to Eqn. (7), where the desired quantity, is the discontinuity of an otherwise real-analytic function. This also takes us back to Eqns. (15a,b): Equation (36) lacks a weight function $\rho(x)$, thus to calculate $\mathcal{F}(z)$ by an N-point Gaussian quadrature requires use of the equivalent quadrature introduced there:

$$\mathcal{F}_N(z) = \sum_{i=1}^N \frac{|<\varphi(x)t(x),\psi(E[i])>|^2}{z-E[i]} w^{EQ}[i]. \tag{37}$$

Note that to calculate $w^{EQ}[i]$ we only need the $E[i]$, assumed to be smoothly interpolated to give $E'[i]$. Thus, there is no need to know anything about a possible OP system underlying this approximation.

But how can this be of utility as, after all, the desired quantities $|<\varphi(x)t(x), \psi(E[i])>|^2$ are still needed? As long as $\varphi(x)t(x)$ is an $L^2$ function, we can imagine that $\psi(E[i])$ only need be known in the support of $\varphi(x)t(x)$, so it can be expanded in a suitable $L^2$ basis set. A convenient way to accomplish this is to replace $\mathcal{H}$, by its matrix representation, $\overline{H}$, in an N-dimensional subset of a complete discrete $L^2$ basis, $\theta_n(r)$. Assuming that diagonalization of gives the N eigenvalues $E[i]$, and $L^2$ eigenfunctions $\psi_{[i]}$ : namely

$$\overline{H}\ \psi_{[i]} = E[i]\ \psi_{[i]}, \text{ with } <\psi_{[i]}, \psi_{[j]}> = \delta_{i,j} \tag{38}$$

The spectral resolution of Eqn. (C) is now replaced by

$$(z - \overline{H})^{-1} = \sum_{i=1}^{N} \psi[i] \frac{1}{(z-E[i])} \psi[i] \tag{39}$$

The analog of Eqn. (E) is now:

$$\overline{F}_N(z) = \sum_{i=1}^{N} \frac{|<\varphi(x)t(x),\psi[i]>|^2}{z-E[i]} . \tag{40}$$

Equating the residues of (E) and (H) gives

$$|<\varphi(x)t(x), \psi(E[i])>|^2 = |<\varphi(x)t(x), \psi[i]>|^2 / w^{EQ}[i], \tag{41}$$

where $w^{EQ}[i] = E'[i]$, and suggesting (correctly) that $w^{EQ}[i]$ is simply, over the support of $\varphi(x)t(x)$, the difference in normalization between the $L^2$ normed $\psi_{[i]}$ and the Dirac $\delta$-function normed $\psi(E[i])$. Proofs of these ideas are the content of Refs. (18) and (36). A worked example now follows.

### B) $L^2$ Discretization in the Absence of a Known OP System: The H-atom Photo-Effect and the Derivative Rule

Figure 11, below, compares the exact Hydrogen atom photo effect cross section, see Refs. (38, 39),

$$\sigma(E) = 2^8/(1+k^2)^5 \text{Exp}[-4\text{ArcTan}(k)/k] / (1-\text{Exp}[-2\pi/k]), \tag{42}$$

to that approximated by discretization of $\mathcal{H}$ and use of the derivative rule. A factor of 1/3, from angular integrations, has been omitted from Eqn. (42) it being an unnecessary part of the present discussion. Here $k = \sqrt{2E}$ is the outgoing momentum (or wave number, in atomic units) of the asymptotic electron of kinetic energy E.

This is the cross section for the process whereby a photon of energy hν ionizes a ground state Hydrogen atom, yielding an asymptotically outgoing electron, e⁻, with asymptotic kinetic energy E = hν – ½, ½ being the binding energy in atomic units. Namely, for the process

$$h\nu + H \rightarrow H^+ + e^-. \qquad (43)$$

This exact cross section is compared with the approximate one making use of eigenfunctions obtained from a discretization of the p-wave, $l = 1$, hydrogenic (Coulomb) Schrödinger operator

$$\mathcal{H} = -\tfrac{1}{2}\, d^2/dr^2 + l(l+1)/2r^2 - 1/r, \qquad (44)$$

discretized in the orthogonal basis of the form, with N = 35, and $\lambda$ a variational parameter:

$$\theta_n(\lambda r) = \mathrm{Exp}[-\lambda r/2]\, (\lambda r)^{l+1}\, L_n^{(2l+2)}(\lambda r), \quad n = 0,1,2,\ldots N-1 \qquad (45)$$

for $l = 1$, and where the $L_n^{2l+2}(\lambda r)$ are the degree n Laguerre polynomials, (Ref. 1b Chapter 18), orthogonal on $r \in [0,\infty)$, with weight function $\mathrm{Exp}[-\lambda r](\lambda r)^{2l+2}$ and measure dr. See Ref. (1b, CH 18, Table 18.3.1). $\varphi(r)t(x)$ is taken to be the exact Hydrogen electronic ground state wave function times the appropriate dipole operator, $t(x)$. See Refs. (38) and (39) for complete discussions of the traditional analytic calculation of this same cross section.

The scaling parameter $\lambda$ was set to 5/2 for the results shown, with N = 35. Varying $\lambda$ from 1/2 to 3 was found to make minimal difference in the quality of results obtained in Fig. 11, except for the fact that they would be obtained for a different set of E[i], as the energy eigenvalues do depend on $\lambda$. The $w^{EQ}[i]$ necessary to extract the cross sections, via Eqn. (41), were obtained by numerical interpolation and differentiation of E[i] as in the Stieltjes inversions of Section 4. Any (reasonable) set of these, with appropriate derivative rule interpolation, would give an essentially identical photo-effect cross section as a function of E continuous on the interval shown, typically to about ten decimal digits

Note that the above choice of basis, Eqn. (45), when used to discretize $\mathcal{H}$ of Eqn. (44) does not give a tri-diagonal J-matrix, nor a three-term recursion allowing it to be mapped onto any known OP system. Thus both the w[i] and $\rho(x)$ are unknown, yet the derivative rule allows computation of the equivalent quadrature weights, and thus approximations to cross-sections. This is the essence of Heller's original ideas, see Refs. (17, 18, 36), where the term *equivalent quadrature* was intended to imply that choice of an L² basis was *equivalent* to introduction of a specific (even if unknown) quadrature type approximation. See Ref. (18) for a similar discretization, with a different, although related, basis which does yield a Pollaczek OP system for suitable $\lambda$. This identification allows an immediate calculation of $w^{EQ}[i]$ in terms of the known weights and weight function of the Attractive Coulomb Pollaczek OP system, but this identification is then subject to the Favard inequality of the discussion following Eqn. (28), and in Ref. (18) itself. Use of the basis of Eqn. (45) is not subject to this restriction.

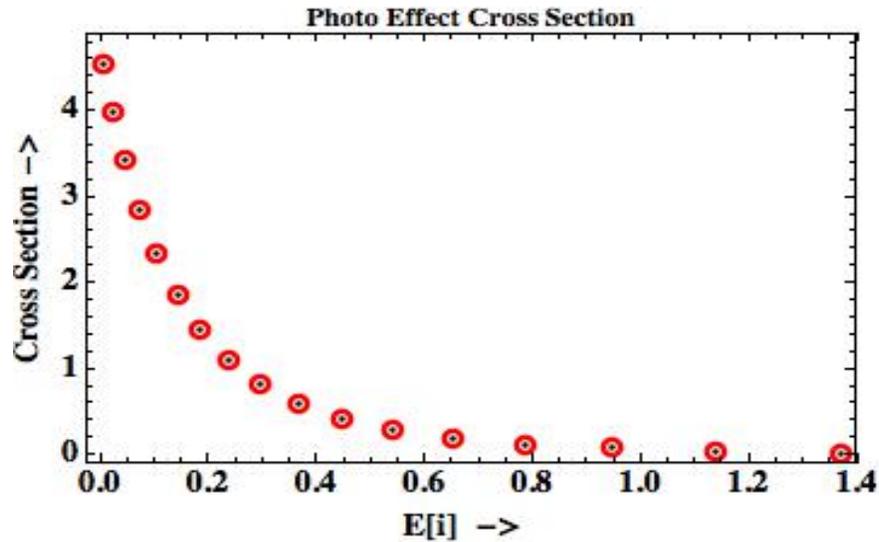

**Figure 11.**
**Caption for Figure 11.** *The photo effect cross section $\sigma(E)$, in atomic units, for single photon ionization of the ground state H atom, as calculated via a discretization in an $L^2$ basis, with N = 35, of the continuum as described above. This is compared, at the eigen-energies E[i] of the discretized Schrödinger operator, to the exact analytic result of Eqn. (42). The black dots are the exact results, and the centers of the red circles are those of the $L^2$ discretization and application of the derivative rule, Eqn. (41), at the same E[n].*

## 6. Summary

The derivative rule, x'[k] = w[k]/ρ(x[k]), is of great utility in Stieltjes inversion and in use with discretized Schrödinger operators. It is exact and universal for the four types of Chebyshev polynomials and validated numerically for a wide range of classical and non-classical OP systems on both finite and infinite intervals. A proof for Nevai-Blumenthal (N-B) class OPs with an a.c. spectrum on [-1,1], has been outlined, and revealed an intriguing restatement of a result known to pure mathematicians: namely that w[k,n]/ρ(x[k,n]) is a universal quantity as $n \to \infty$ for N-B OP systems, of possible interest to applied and numerically oriented users of mathematics. How this latter extends most effectively to OP systems, not of the N-B class, such as the classical Laguerre OPs, or to a mixture of *a.c.* and *discrete* spectra are open questions.

## 7. Acknowledgements

My former students, i.e. my teachers, from the 1970s, David Oxtoby, Tom Rescigno, Terry Murtaugh, Lou Fishman, John Broad, S. V. ONeil, and especially Hashim Yamani and Eric J. Heller are most gratefully acknowledged.  Rather more recently, going on 45 years later, I have benefited from conversations with Nico Temme, Barry Schneider, Tom Koornwinder, Adri Olde Daalhuis, and Mourad Ismail.  These being members of the groups revising and updating the NIST Digital Library of Mathematical Functions as fellow members of the DLMF Senior Editorial Board, Chapter Authors, or NIST Staff.   In particular: Tom Koornwinder suggested the Chihara additions to Ref. (24), and Mourad Ismail pointed out several important additions to those originally listed in Ref. (29).  The support of the present work, via a NIST grant for work on development of the Orthogonal Polynomial chapter of the DLMF, is gratefully acknowledged.

## Appendix: The Chebyshev Polynomials, a Brief Overview

The quadrature abscissas, $x[k]$ (= $x[k,n]$), and their derivatives $x'[k,n]$ with respect to "k," weights, and weight functions, as used in the text, are stated for the $n^{th}$ order Chebyshev polynomials of the first through the fourth kinds, labelled 1 through 4. See Ref. (1a,b). Minus signs, not often present, insure that the $x[k,n]$ are increasing functions of "k" making the derivatives positive, see the notations as laid out in **Section 2A**.  Universality of the $x'[k,n]$ is evident as $n \to \infty$.

**First Kind**

$x1[k, n] = -\mathrm{Cos}((2k - 1)*\pi/(2n))$      (A.1.1)

$x1'[k, n] = \pi/n \, \mathrm{Sin}((2k - 1)\pi)/(2n)) = (\pi/n)\sqrt{(1 - x1[i,n]^{\wedge}2)}$      (A.1.2)

$w1[k, n] = \pi/n$      (A.1.3)

$\rho C1(x) = 1/\sqrt{(1 - x^{\wedge}2)}, \quad x \in [-1,1]$      (A.1.4)

**Second Kind**

$x2[k, n] = -\mathrm{Cos}(k\pi/(n + 1))$      (A.2.1)

$x2'[k\,n] = \pi/(1 + n)\,\mathrm{Sin}((k\pi)/(n + 1) = \pi/(1 + n)\sqrt{(1 - x2[k,n]^{\wedge}2)}$      (A.2.2)

$w2[k, n] = \pi/(n + 1)\,\mathrm{Sin}(i\,\pi/(n + 1))^{\wedge}2$      (A.2.3)

$\rho C2(x) = \sqrt{(1 - x^{\wedge}2)}, \, x \in [-1,1]$      (A.2.4)

**Third Kind**

$$x3[k, n] = -\cos(2k\pi/(2n + 1)) \quad \text{(A.3.1)}$$

$$x3'[k, n] = 2\pi/(2n + 1) \sin((2k\pi)/(2n + 1)) = 2\pi/(2n + 1) \sqrt{(1 - x3[k, n]^2)} \quad \text{(A.3.2)}$$

$$w3[k,n] = 4\pi/(2n+1) \sin((n-(k-1))\pi/(2n + 1))^2 \quad \text{(A.3.3)}$$

$$\rho C3(x) = \frac{\sqrt{(1-x)}}{\sqrt{(1+x)}} \quad , x \in [-1,1] \quad \text{(A.3.4)}$$

**Fourth Kind**

$$x4[k, n] = -\cos(2k\pi/(2n + 1)) \quad \text{(A.4.1)}$$

$$x4'[k, n] = 2\pi/(2n + 1) \sin((2k\pi)/(2n + 1)) = 2\pi/(2n + 1) \sqrt{(1 - x4[k, n]^2)} \quad \text{(A.4.2)}$$

$$w4[k,n] = 4\pi/(2n+1) \sin(k\pi/(2n + 1))^2 \quad \text{(A.4.3)}$$

$$\rho C4(x) = \frac{\sqrt{(1 \mp x)}}{\sqrt{(1 - x)}} \quad , x \in [-1,1] \quad \text{(A.4.4)}$$

**Universality of the x'[k,n]**

Note that as $n \to \infty$, and for k not too close to either 1 or n:

$$x1'[k,n] = x2'[k,n] = x3'[k,n] = x4'[k,n]$$

$$= (\pi/n) \sqrt{(1 - \cos(k\pi/n)^2)} = (\pi/n) \sqrt{(1 - x[k, n]^2)}, \quad \text{(A.5)}$$

where $x[k,n] = -\cos(k\pi/n)$.

never thought of the derivative rule as being anything other than a density of states approximation to be used to calculate the equivalent quadrature weights (see Refs. 18, and 36, below) for use in approximate, discrete, spectral resolutions of Schrödinger resolvents of the type as further discussed in Section 5.  The quadrature abscissas in those discussions were energies, and the # of quantum states, di, per energy interval dE, being a quantum state density, thus di/dE, or $(dE/di)^{-1}$.  Heller's original idea is also related to similar ideas, then current in theoretical atomic physics, as described in U. Fano and J. W. Cooper, *Spectral Distributions of Atomic Oscillator Strengths,* Rev. Mod. Phys. **40**,411(1968), in following oscillator strength distributions for highly excited quantum bound states (i.e. Rydberg states) merging smoothly into the Coulomb photo ionization (or photo effect) continuum. See also: *Atomic Collisions and Spectra*, U. Fano and A. R. P. Rau, Academic Press, New York, 1986, pp 33-35, where a related use of a derivative is introduced.